\documentstyle[12pt]{article}
\begin{document}
\newpage
\pagestyle{empty}
\setcounter{page}{0}
%

\renewcommand{\theequation}{\thesection.\arabic{equation}}
\newcommand{\sect}[1]{\setcounter{equation}{0}\section{#1}}

\newfont{\twelvemsb}{msbm10 scaled\magstep1}
\newfont{\eightmsb}{msbm8}
\newfont{\sixmsb}{msbm6}
\newfam\msbfam
\textfont\msbfam=\twelvemsb
\scriptfont\msbfam=\eightmsb
\scriptscriptfont\msbfam=\sixmsb
\catcode`\@=11
\def\Bbb{\ifmmode\let\next\Bbb@\else
  \def\next{\errmessage{Use \string\Bbb\space only in math mode}}\fi\next}
\def\Bbb@#1{{\Bbb@@{#1}}}
\def\Bbb@@#1{\fam\msbfam#1}
\newfont{\twelvegoth}{eufm10 scaled\magstep1}
\newfont{\tengoth}{eufm10}
\newfont{\eightgoth}{eufm8}
\newfont{\sixgoth}{eufm6}
\newfam\gothfam
\textfont\gothfam=\twelvegoth
\scriptfont\gothfam=\eightgoth
\scriptscriptfont\gothfam=\sixgoth
\def\frak{\frak@}
\def\frak@#1{{\fam\gothfam{{#1}}}}
\def\frak@@#1{\fam\gothfam#1}
\catcode`@=12
%
%
%
\def\CC{{\Bbb C}}
\def\NN{{\Bbb N}}
\def\QQ{{\Bbb Q}}
\def\RR{{\Bbb R}}
\def\ZZ{{\Bbb Z}}
\def\cA{{\cal A}}          \def\cB{{\cal B}}          \def\cC{{\cal C}}
\def\cD{{\cal D}}          \def\cE{{\cal E}}          \def\cF{{\cal F}}
\def\cG{{\cal G}}          \def\cH{{\cal H}}          \def\cI{{\cal I}}
\def\cJ{{\cal J}}          \def\cK{{\cal K}}          \def\cL{{\cal L}} 
\def\cM{{\cal M}}          \def\cN{{\cal N}}          \def\cO{{\cal O}}
\def\cP{{\cal P}}          \def\cQ{{\cal Q}}          \def\cR{{\cal R}} 
\def\cS{{\cal S}}          \def\cT{{\cal T}}          \def\cU{{\cal U}}
\def\cV{{\cal V}}          \def\cW{{\cal W}}          \def\cX{{\cal X}}
\def\cY{{\cal Y}}          \def\cZ{{\cal Z}}
\def\qed{\hfill \rule{5pt}{5pt}}
\def\id{\mbox{id}}
\def\ggo{{\frak g}_{\bar 0}}
\def\uqggo{\cU_q({\frak g}_{\bar 0})}
\def\uqggp{\cU_q({\frak g}_+)}
\def\typeA{{\em type $\cA$}}
\def\typeB{{\em type $\cB$}}
\newtheorem{lemma}{Lemma}
\newtheorem{prop}{Proposition}
\newtheorem{theo}{Theorem}

$$
\;
$$
\rightline{CPTH-S458.0696}
\rightline{q-alg/9610031}
\rightline{June 96}

\vfill
\vfill
\begin{center}

  {\LARGE {\bf {\sf On ${\cal U}_h(sl(2))$, ${\cal U}_{h}(e(3))$ and their
Representations}}}
 \\[2cm]

\smallskip 

{\large B. Abdesselam\footnote{abdess@orphee.polytechnique.fr},
A. Chakrabarti\footnote{chakra@orphee.polytechnique.fr} and 
R. Chakrabarti\footnote{Permanent address: Department of Theoretical Physics, University of Madras, Guindy Campus, Madras-600025, India}}

\smallskip 

\smallskip 

\smallskip 

{\em  \footnote{Laboratoire Propre du CNRS UPR A.0014}Centre de Physique 
Th\'eorique, Ecole Polytechnique, \\
91128 Palaiseau Cedex, France.}

\end{center}

\vfill

\begin{abstract}
By solving a set of recursion relations for the matrix elements of the 
${\cal U}_h(sl(2))$ generators, the finite dimensional highest weight 
representations of the algebra were obtained as factor representations. 
Taking a nonlinear combination of the generators of the two copies of the 
${\cal U}_h(sl(2))$ algebra, we obtained ${\cal U}_h(so(4))$ algebra.    
The latter, on contraction, yields ${\cal U}_h(e(3))$ algebra. A nonlinear 
map of ${\cal U}_h(e(3))$ algebra on its classical analogue $e(3)$ 
was obtained. The inverse mapping was found to be singular. It signifies 
a physically interesting situation, where in the momentum basis, a 
restricted domain of the eigenvalues of the classical operators is mapped 
on the whole real domain of the eigenvalues of the deformed operators.
   
\end{abstract}

\vfill
\vfill

\newpage
\pagestyle{plain}

\sect{Introduction} 

The enveloping Lie algebra ${\cal U}(sl(2))$ has two distinct quantizations: 
The first one is called the Drinfeld-Jimbo deformation (standard 
$q$-deformation) \cite{Drin,Jim}, whereas the second one is called the 
Jordanian deformation (nonstandard $h$-deformation) \cite{Manin,Ohn} and 
may be obtained as a contraction of the Drinfeld-Jimbo one \cite{Iran}. 
Recently there is much interest in studies relating to various aspects of 
the $h$-deformed algebra ${\cal U}_{h}(sl(2))$. In particular, a two 
parametric deformation of the dual algebra ${\cal F}un_{h,h'}(GL(2))$ 
was obtained in \cite{Agha}. This author also constructed \cite{Agha} 
the differential calculus in the quantum plane. Quantum de Rham complexes 
associated with the $h$-deformed algebra ${\cal F}un_{h}(sl(2))$ was given in 
\cite{Kiri}. The universal ${\cal R}$-matrix of the algebra ${\cal U}_{h}
(sl(2))$ was obtained \cite{Vlad,Shar}. Various non-semisimple $h$-deformed 
algebras were constructed at contraction limits \cite{Ball,Ball2,Shar,para}. 
The $h$-deformation was also extended to the case of superalgebras 
\cite{Dab}.   

One of the studies made in the present article is as follows. Using the 
standard singular vector construction method, we study the finite dimensional 
highest weight representations of ${\cal U}_h(sl(2))$ algebra. Along similar 
lines, investigations were made before in \cite{Dob}. Compared to 
\cite{Dob}, a distinctive feature in our approach is that we develop a set of 
recursion relations, which may be easily solved to determine the matrix
elements of the generators of the algebra. These matrix elements, in turn,
specify the singular vectors leading to finite dimensional irreducible 
representations. This representations may also be obtained by exploiting 
a recently \cite{ACC} proposed nonlinear invertible map between the 
generators of
${\cal U}_{h}(sl(2))$ and the classical $sl(2)$ generators. We wish to
stress, however, that the construction of the finite dimensional highest 
weight representations via the singular vector technique relies on first 
principles; and, therefore, may be useful for other nonlinear algebras, 
which exhibit no such maps on the corresponding linear algebras. In a 
continuation of our studies on the ${\cal U}_{h}(sl(2))$ algebra \cite{ACC}, 
we consider the contracted ${\cal U}_{h}(e(2))$ algebra \cite{Ball2}. The 
latter may be also be mapped on the classical algebra $e(2)$.                 

We also consider the ${\cal U}_h(e(3))$ algebra obtained as a contraction of 
the ${\cal U}_h(so(4))$ algebra \cite{Ball,Shar}, that may be realized from 
two copies of ${\cal U}_{h}(sl(2))$. It is of interest that, parallel to the
scenario discussed for ${\cal U}_{h}(sl(2))$ algebra, ${\cal U}_h(e(3))$
algebra may also be realized in terms of the $e(3)$ generators. This 
enormously simplifies the problem of finding the irreducible representations 
of ${\cal U}_h(e(3))$ algebra . In contrast to a similar map obtained  
in \cite{Ac4} for ${\cal U}_q(e(3))$ algebra and also the realization 
\cite{Kosi} of the $\kappa$-Poincar\'e algebra \cite{Luk1,Luk2} in terms of
the classical Poincar\'e generators, the corresponding map for the present 
${\cal U}_h(e(3))$ algebra on classical $e(3)$ exhibits a singularity. This 
may be of physical interest.

Let $h$ be an arbitrary complex parameter. The 
algebra ${\cal U}_{h}(sl(2))$ is then an associative algebra over $\CC$
generated by $H$, $X$ and $Y$, satisfying the commutation relations 
\cite{Ohn}      
\begin{eqnarray}
  && [H,X]= 2 {\sinh hX \over h},\nonumber \\
  && [H,Y]=-Y(\cosh hX)-(\cosh hX)Y,\nonumber\\  
  && [X,Y]=H. 
\end{eqnarray}
The coalgebra structure of  ${\cal U}_{h}(sl(2))$ reads \cite{Ohn}
\begin{eqnarray}
&& \bigtriangleup (X)=X \otimes 1 + 1 \otimes X,\nonumber  \\
&& \bigtriangleup (Y)=Y\otimes e^{h X} + e^{-h X} 
 \otimes Y, \nonumber\\
&& \bigtriangleup (H)=H \otimes e^{h X} + e^{-h X} 
 \otimes H, \nonumber\\
&& \varepsilon (X)= \varepsilon (Y)= \varepsilon (H)=0,\nonumber \\
&& S(X) =- X,\qquad
 S(Y) =- e^{h X}Y e^{-h X},\qquad
S(H) =- e^{h X}H e^{-h X}.  
\end{eqnarray}                    
The Casimir element of ${\cal U}_{h}(sl(2))$ is given by \cite{Ball}  
\begin{eqnarray}
&&C={1\over 2h}\biggl(Y(\sinh hX)+(\sinh hX)Y\biggr)+{1\over 4}H^2+
{1\over 4}(\sinh hX)^2.
\end{eqnarray}

\sect{Representations of ${\cal U}_h(sl(2))$ algebra} 

The finite dimensional highest 
weight representations of ${\cal U}_h(sl(2))$ as the factor-representations 
of the corresponding Verma modules were considered in \cite{Dob}. The 
factorization $\;$scheme was $\;$carried out by $\;$the standard 
$\;$singular vector 
construction method. In \cite{Dob}, the operators $H$, $Y$, $\cosh (hX)$ 
and $\sinh(hX)$ were chosen as generators. We follow the same route here. 
However, a 
different choice of the generators of the algebra ${\cal U}_h(sl(2))$ 
allows us to express the matrix elements in terms of a simple set of recursion 
relations which may be 
easily solved. These matrix elements, in turn, completely determine the 
singular vectors at arbitrary levels. The appearance of a singular vector
in a Verma module signals its reducibility. To
obtain a finite dimensional irreducible representation, the submodule 
generated by treating the singular vector as a highest weight vector must be 
factored off. This scheme provides all irreducible representations at 
arbitrary dimensions. Alternately, these irreducible representations may be 
obtained via a recently developed \cite{ACC} nonlinear invertible map 
between the ${\cal U}_h(sl(2))$ generators and the generators of the classical 
undeformed ($h=0$) $sl(2)$ algebra. We will briefly elucidate this procedure 
later.

The highest weight vector $w_{0}$, where $\lambda$ is the highest weight, 
satisfies the relations
\begin{eqnarray}
X.w_{0}=0,\qquad\qquad H.w_{0}=\lambda w_{0}.
\end{eqnarray}     
The Verma module $M$ is generated by the repeated actions of $Y$ on $w_{0}$
as
\begin{eqnarray}
w_{n}=Y^{n}.w_{0},\qquad\qquad n\in \NN.
\end{eqnarray}             
Using the commutations relations $(1.1)$, it is evident that the actions of 
$X$ and $H$ on the vector space are described as
\begin{eqnarray}
&& X.w_{n}=\sum_{k=0}^{[(n-1)/2]}X_{n}^{n-1-2k}\; w_{n-1-2k},  \\
&& H.w_{n}= \sum_{k=0}^{[n/2]}H_{n}^{n-2k}\; w_{n-2k}, 
\end{eqnarray}       
where $H^{0}_{0}=\lambda$ and $[x]$ denotes the integer part in $x$. We 
note here that the actions of $X$ and $H$ on a
state in the Verma module create the sequences of states differing in their 
indices by two. Our task is now to develop the recursion relations between 
the above matrix elements. To this end, we use $(2.2)$ to obtain 
\begin{eqnarray}
&& X.w_{n}=([X,Y]+YX)w_{n-1},  \\
&& H.w_{n}=([H,Y]+YH)w_{n-1}. 
\end{eqnarray}        
The commutation relations $(1.1)$ may now be exploited to obtain the following 
recursion relations
\begin{eqnarray}
&& X_{n}^{m}=X_{n-1}^{m-1}+H_{n-1}^{m},  \\
&& H_{n}^{n}= H_{n-1}^{n-1}-2,  \\
&& H_{m+2n}^{m}= H_{m+2n-1}^{m-1}-
\sum_{k=1}^{n}{h^{2k}\over (2k)!} \sum_{\delta=0}^{1}\sum_{\{\Delta_{i}
|i=(1,2,\cdots,2k)\}} ' Z_{m,2n,\delta,\{\Delta_{i}\}}, \qquad\qquad 
\qquad \qquad 
\end{eqnarray} 
where the primed sum in the rhs is performed over the following all possible 
partitions of $2n$ among an even number of positive odd integers 
\begin{eqnarray}
\Delta_1+\Delta_2+\cdots + \Delta_{2k}=2n,\qquad \Delta_{i}\;mod\;2=1\;\;
\hbox{for}\;\; i=(1,2,\cdots, 2k)
\end{eqnarray}
and
\begin{eqnarray}
Z_{m,2n,\delta,\{\Delta_{i}\}}= X_{m+2n-\delta}^{m+2n-\delta-\Delta_{1}}
 X_{m+2n-\delta-\Delta_1}^{m+2n-\delta-\Delta_{1}-\Delta_{2}}\cdots
X^{m-\delta}_{m+2n-\delta-\Delta_{1}-\Delta_{2}-\cdots - \Delta_{2k-1}}.
\end{eqnarray}
The partition $(2.10)$ reduces the problem of finding the matrix elements 
of $H$ and 
$X$ to a combinatorial exercise. This permits us to write down, starting from
the known values at ($h=0$) limit, the matrix elements of these operators at 
arbitrary dimensions. We demonstrate this in the present section. To cite an 
example
of the partition $(2.10)$, we enlist all possible cases of partitioning $8$ 
($=2n$): $\{(7,1),(1,7)$; $(5,3),(3,5)$; $(5,1,1,1)$ and permutations; 
$(3,3,1,1)$ and permutations; $(3,1,1,1,1,1)$ and permutations; 
$(1,1,1,1,1,1,1,1)\}$. Rewriting $(2.7)$ differently, we obtain
\begin{eqnarray}
X_{m+2n+1}^{m}=\sum_{k=0}^{m}H_{k+2n}^{k}, \qquad\qquad n\geq 0.
\end{eqnarray}
The matrix elements $\sim O(1)$ surviving in the classical limit 
($h = 0$) can be immediately obtained from $(2.8)$ and $(2.12)$
\begin{eqnarray}
&& H_{n}^{n}=\lambda-2n ,  \\
&& X_{n+1}^{n}=(n+1)(\lambda-n). 
\end{eqnarray}
From (2.9), (2.12), (2.13) and (2.14), it is evident that the matrix elements 
$X_{m+2n+1}^{m}$ and $H_{m+2n}^{m}$ are $\sim O(h^{2n})$. We now demonstrate 
that starting with the known matrix elements $(2.13)$ and $(2.14)$ at $\sim
O(1)$, we may explicitly evaluate the matrix elements appearing at an 
arbitrary power of $h$ by solving the recursion relations (2.9) and (2.12). 
To determine the matrix elements $H_{n+2}^{n}\sim O(h^2)$, we obtain from 
(2.9),
\begin{eqnarray}
H_{n+2}^{n}=-h^2\sum_{k=0}^{n-1}X_{k+2}^{k+1}X_{k+1}^{k} -
{h^2\over 2} X_{n+2}^{n+1}X_{n+1}^{n}.
\end{eqnarray}
Using the explicit values $(2.14)$ of the elements $X_{n+1}^{n}$, we obtain 
\begin{eqnarray}
H_{n+2}^{n}= h^2 \rho^{n}_{(2)}, 
\end{eqnarray}             
where
\begin{eqnarray}
\rho^{n}_{(2)}= -\sum_{k=0}^{n-1}(k+1)(k+2)(\lambda-k)(\lambda-k-1)
-{1\over  2} (n+1)(n+2)(\lambda-n)(\lambda-n-1).
\end{eqnarray} 
The elements $X_{n+3}^{n} \sim O(h^2)$ now follow from (2.12) and (2.16) as
\begin{eqnarray}
X_{n+3}^{n} =h^2\sigma_{(2)}^{n},  
\end{eqnarray} 
where
 \begin{eqnarray}
\sigma_{(2)}^{n}=\sum_{k=0}^n \rho_{(2)}^{k}.
 \end{eqnarray} 
This completes our explicit evaluation of all matrix elements $\sim O(h^2)$. 
Exploiting these explicitly known elements ($\sim O(1)$, $O(h^2)$) in (2.13), 
(2.14), (2.16) and (2.18), we now determine the elements $H_{n+4}^{n}$ and
$X_{n+5}^{n}\sim O(h^4)$. We enlist the result as follows
\begin{eqnarray}
H_{n+4}^{n}=h^4 \rho_{(4)}^{n}, 
\end{eqnarray}          
where 
\begin{eqnarray}
&&\rho^{n}_{(4)}= -\sum_{k=0}^{n-1}\biggl( (k+4)(\lambda -k-3) 
\sigma_{(2)}^{k}+ (k+1)(\lambda-k) \sigma_{(2)}^{k+1}
+ 2 (4!) \biggl({k+4 \atop 4}\biggr)\biggl({\lambda -k \atop \lambda-k 
-4}\biggr) \biggr)\nonumber \\
&&\phantom{\rho^{n}_{(4)}= }  -{1\over 2} \biggl( (n+4)(\lambda -n-3) 
\sigma_{(2)}^{n}+ (n+1)(\lambda-n) \sigma_{(2)}^{n+1}
+ 2 (4!) \biggl({n+4 \atop 4}\biggr)\biggl({\lambda -n \atop \lambda-n 
-4}\biggr) \biggr) \nonumber \\
&& 
\end{eqnarray} 
and
\begin{eqnarray}
X_{n+5}^{n} =h^4\sigma_{(4)}^{n},  
\end{eqnarray} 
where
 \begin{eqnarray}
\sigma_{(4)}^{n}=\sum_{k=0}^n \rho_{(4)}^{k}. 
 \end{eqnarray} 
The general scheme of determination of an arbitrary matrix element is 
evident now. Assuming that all the matrix elements appearing up to the order 
$\sim O(h^{2(n-1)})$ have been explicitly determined, we proceed to determine 
the matrix elements $H_{m+2n}^{m}$ and $X_{m+2n+1}^{m}$ ($\sim O(h^{2n})$). 
The relation (2.9) may be reorganized as follows
\begin{eqnarray}
H_{m+2n}^{m}=-\sum_{k=1}^{n} {h^{2k} \over (2k)!} 
\sum_{\{\Delta_i| i=(1,\cdots,2k)\}}' \biggl(2\sum_{l=0}^{m-1}
Z_{l,2n,0,\{\Delta_i\}}+Z_{m,2n,0,\{\Delta_{i}\}}\biggr), 
\end{eqnarray}
where the primed sum in the rhs has been described in (2.10). The matrix 
elements appearing in the rhs of (2.24) are already known and may be utilized 
to explicitly evaluate $H_{m+2n}^{m}$. Equations $(2.12)$ and $(2.24)$ now 
determine $X_{m+2n+1}^{m}$. This completes our determination of all 
the matrix elements of the operator $X$ and $H$ acting on the Verma module.

For $\lambda=2j$ ($j=0,{1\over 2},1 ,\cdots$), singular vectors $\{w_{s}^{
(j)}\;|\; (2j+1)\in \NN\}$ appear in the Verma module, which is, thus, 
rendered reducible. An irreducible representation of $(2j+1)$ dimension 
is obtained by taking the quotient module $L^{j}=M/\{{\cal U}_h(sl(2)).
w_{s}^{(j)}\}$. The singular vectors $\{w_{s}^{(j)}\}$ are annihilated by 
$X$ and are eigenvectors of $H$ 
\begin{eqnarray}
&& X.w_{s}^{(j)}=0,\\
&& H.w_{s}^{(j)}=\lambda^{(j)} w_{s}^{(j)}.    
\end{eqnarray}
We note that $(2.25)$, together with the commutation relation (1.1), demands 
(2.26) to be satisfied. The zero mode condition (2.25) gives rise to a set of 
linear equations determining the singular vectors $w_{s}^{(j)}$ completely.
We describe this bellow. The  matrix elements obtained before suggest an 
ansatz for $\{w_{s}^{(j)}\}$
\begin{eqnarray}
w_{s}^{(j)}=w_{2j+1}+\sum_{p=1}^{[j]} C_{p}^{(j)}w_{2j-2p+1}.
\end{eqnarray} 
Substiting (2.27) in (2.25), we get
\begin{eqnarray}
(2j+1)(\lambda-2j)w_{2j} +\sum_{r=1}^{[j]} X_{2j+1}^{2j-2r}w_{2j-2r}
+\sum_{p=1}^{[j]} C_{p}^{(j)}\sum_{r=0}^{[j]-p} X_{2j+1-2p}^{2j-2p-2r}
w_{2j-2r-2p}=0.
\end{eqnarray} 
For $\lambda=2j$, the first term in $(2.28)$ vanishes. Making a change in the 
index in the second summation, we obtain
\begin{eqnarray}
\sum_{r=1}^{[j]} X_{2j+1}^{2j-2r}w_{2j-2r}
+\sum_{p=1}^{[j]} C_{p}^{(j)}\sum_{r=p}^{[j]} X_{2j+1-2p}^{2j-2r}
w_{2j-2r}=0.
\end{eqnarray} 
In (2.29), it is understood that the matrix elements of $X$ are evaluated at 
$\lambda=2j$. Reversing the order of summation in the second term (2.29), 
we get 
\begin{eqnarray}
\sum_{r=1}^{[j]}\biggl(X_{2j+1}^{2j-2r}+\sum_{p=1}^{r}C_{p}^{(j)}
X_{2j+1-2p}^{2j-2r}\biggr)w_{2j-2r}=0.
\end{eqnarray} 
This enforces the coefficients $C_{p}^{(j)}$ ($p=1,\cdots,[j]$)
to satisfy $[j]$ linear (in fact, triangular) equations, where the 
matrix elements of $X$ are evaluated at $\lambda=2j$:
\begin{eqnarray}
X_{2j+1}^{2j-2r}(\lambda=2j)+\sum_{p=1}^{r}C_{p}^{(j)}
X_{2j+1-2p}^{2j-2r} (\lambda=2j) =0, \qquad r=(1,\cdots, [j]).
\end{eqnarray}         
The equations (2.31) may be easily solved. We complete our discussion of the 
singular vectors corresponding to arbitrary $(2j+1)$ dimensional 
representations by enlisting the first few singular vectors for different 
choices of $\lambda$ 
\begin{equation}
\begin{array}{ll}
\hbox{\sl Value of $\lambda$ ($=2j$)}  &\qquad  \hbox{\sl Singular vector} 
\qquad\qquad \\
\qquad 0 & \qquad w_1 \\
\qquad 1  & \qquad w_2 \\
\qquad 2 & \qquad w_3+h^2 w_1 \\
\qquad 3 & \qquad w_4+ 6h^2 w_2 \\
\qquad 4 & \qquad w_5+ 21 h^2 w_3+ 36 h^4 w_1 \\
\qquad 5 & \qquad w_6+ 56 h^2w_4 +460 h^4 w_2 \\
\qquad 6 & \qquad w_7+ 126h^2 w_5 +3105 h^4 w_3 + 8100 h^6 w_1 \\
\qquad 7 & \qquad w_8 +252 h^2 w_6 + 14796 h^4 w_4 + 166320 h^6 w_2 
\end{array}
\end{equation}   
The analytical expression derived in \cite{Dob} allows one to find the 
singular vectors for $\lambda \leq 3$. Using REDUCE, this author also 
derived expressions necessary for finding singular vectors up to $\lambda 
\leq 6$. In our recursive scheme, the problem of determination of singular 
vectors is much
simplified. Determination of the matrix elements of $X$ by our recursive 
method also yields automatically the singular vectors at arbitrary levels. 
The matrix element of the operator $Y$ may, now, be readily determined. The 
recipe is as follows. In order to extract the irreducible representations 
of dimension $(2j+1)$, the singular vector $w_{s}^{(j)}$ existing for 
$\lambda=2j$ is identified with the null vector
\begin{eqnarray} 
w_{s}^{(j)}\approx 0.
\end{eqnarray}  
This, in conjunction with (2.2) and (2.27), now leads to  
\begin{eqnarray} 
Y.w_{2j}=-\sum_{p=1}^{[j]}C_{p}^{(j)}w_{2j-2p+1}.  
\end{eqnarray} 
Our construction of states (2.2) and the coefficients $\{C_{p}^{(j)}\}$, 
already determined by the system of linear equations (2.31), now provide 
the all the matrix elements of the operator $Y$. 

We supplement our preceding description with an example, where we construct 
the irreducible representations for $j={7\over 2}$ case. The factor module 
now consists of the vector space $\{w_{i}|i=0,\cdots,7\}$, which may be 
identified as follows
\begin{eqnarray}
(w_{i})_{j}=\delta_{ij}, \qquad\qquad (i,j)=(0,1,\cdots, 7).  
\end{eqnarray}
Using our construction of the matrix elements (2.12), (2.13), (2.14) 
and (2.24) we obtain the representations of the operators $X$ and $H$. 
The relations (2.2), (2.31) and (2.34) provide the representations of 
the operator $Y$. The $j=7/2$ representation reads: 
\begin{eqnarray}
&& X=\pmatrix{
0 & 7 & 0 & -42h^2 & 0 & 252 h^4 & 0 &  58968 h^6 \cr
0 & 0 & 12 & 0 & -216 h^2  & 0 & 4176 h^4 & 0 \cr
0 & 0 & 0 & 15 & 0 & -600 h^2 & 0 & 22500 h^4 \cr
0 & 0 & 0 & 0 & 16 & 0 & -1224h^2 & 0 \cr
 0 & 0 & 0 & 0 & 0 & 15 & 0 & -2058h^2 \cr
 0 & 0 & 0 & 0 & 0 & 0 & 12 & 0 \cr
   0 & 0 & 0 & 0 & 0 & 0 & 0 & 7 \cr
   0 & 0 & 0 & 0 & 0 & 0 & 0 & 0\cr}, \nonumber \\
&& Y=\pmatrix{
0 & 0 & 0 & 0 & 0 & 0 & 0 &  0 \cr
1 & 0 & 0 & 0 & 0 & 0 & 0 & 0 \cr
0 & 1 & 0 & 0 & 0 & 0 & 0 & -166320h^6 \cr
0 & 0 & 1 & 0 & 0 & 0 & 0 & 0 \cr
0 & 0 & 0 & 1 & 0 & 0 & 0 & -14796h^4 \cr
0 & 0 & 0 & 0 & 1 & 0 & 0 & 0 \cr
0 & 0 & 0 & 0 & 0 & 1 & 0 & -252h^2 \cr
0 & 0 & 0 & 0 & 0 & 0 & 1 & 0\cr}, \nonumber \\
&& H=\pmatrix{
7 & 0 & -42h^2 & 0       & 252h^4    &    0   &   58968h^6 &  0 \cr
0 & 5 & 0      & -174h^2 & 0         & 3924h^4  & 0 & 88776h^6 \cr
0 & 0 & 3      & 0       & -384h^2  & 0      & 18324h^4 & 0 \cr
0 & 0 & 0      & 1       & 0         & -624h^2  & 0 & 49212 h^4\cr
0 & 0 & 0      & 0       & -1        & 0    & -834 h^2 & 0 \cr
0 & 0 & 0      & 0       & 0         & -3      & 0  & -966 h^2 \cr
0 & 0 & 0      & 0       & 0         & 0      & -5   &   0  \cr
0 & 0 & 0      & 0       & 0         & 0      & 0    &  -7  \cr}.
\end{eqnarray}
For later comparison purpose, we also express the $j=7/2$ 
representation of the ${\cal U}_h(sl(2))$ algebra, where the operator $H$ 
has been diagonalized. It reads            
\begin{eqnarray}
&& X=\pmatrix{
0 & 7 & 0 & 105h^2 & 0 & 3780 h^4 & 0 &  56700 h^6 \cr
0 & 0 & 12 & 0 & 240 h^2  & 0 & 6480 h^4 & 0 \cr
0 & 0 & 0 & 15 & 0 & 300 h^2 & 0 & 3780 h^4 \cr
0 & 0 & 0 & 0 & 16 & 0 & 240 h^2 & 0 \cr
 0 & 0 & 0 & 0 & 0 & 15 & 0 & 105 h^2 \cr
 0 & 0 & 0 & 0 & 0 & 0 & 12 & 0 \cr
   0 & 0 & 0 & 0 & 0 & 0 & 0 & 7 \cr
   0 & 0 & 0 & 0 & 0 & 0 & 0 & 0\cr}, \nonumber \\
&& Y=\pmatrix{
0 & \displaystyle -{21 h^2 \over 2} & 0 & \displaystyle {315 h^4 \over 4}   
& 0 & 4725h^6 & 0 &  99225h^8 \cr
1 & 0 & -33h^2 & 0 & 180h^4 & 0 & 8100 h^6 & 0 \cr
0 & 1 & 0 &  \displaystyle -{105 h^2 \over 2}   & 0 & 225h^4 & 0 
& 4725h^6 \cr
0 & 0 & 1 & 0 & -60h^2 & 0 & 180h^4 & 0 \cr
0 & 0 & 0 & 1 & 0 & \displaystyle -{105 h^2 \over 2}  & 0 &
\displaystyle {315 h^4 \over 4}  \cr
0 & 0 & 0 & 0 & 1 & 0 & -33h^2 & 0 \cr
0 & 0 & 0 & 0 & 0 & 1 & 0 & \displaystyle -{21 h^2 \over 2} \cr
0 & 0 & 0 & 0 & 0 & 0 & 1 & 0\cr}, \nonumber \\
&& \nonumber \\
&& H=\hbox{diag}( 
7,\;5,\;3,\;1,\;-1,\;-3,\;-5,\;-7).
\end{eqnarray}
This completes our recipe for constructing an arbitrary $(2j+1)$ dimensional 
irreducible representation of the algebra ${\cal U}_h(sl(2))$ by the process 
of factorization of the Verma module $M$. The above method, developed from 
the first principles, may now be compared with an alternative way of obtaining 
the irreducible representations by exploiting a recently developed \cite{ACC} 
invertible map between the generators of the  ${\cal U}_h(sl(2))$ 
algebra and the undeformed $sl(2)$ generators. The map reads
\begin{eqnarray}   
&& X = {2\over h}\;\hbox{arctanh}({h\;J_{+}\over 2}), \nonumber \\
&& Y=\sqrt{1-{h^2J_{+}^2 \over 4}}\;J_{-}\; \sqrt{1-{h^2J_{+}^2 \over 4}},\nonumber \\
 && H=J_{0},
\end{eqnarray}
where $(J_{\pm}, J_{0}$) satisfy the $sl(2)$ algebra
\begin{eqnarray}
  && [J_{0},J_{\pm}]= \pm 2 J_{\pm}, \qquad\qquad\qquad
   [J_{+},J_{-}]=J_{0}
\end{eqnarray}
and the following cocommutative coproduct relations 
\begin{eqnarray}
\bigtriangleup(J_{i})=J_{i}\otimes 1+1\otimes J_{i}, \qquad i=(\pm ,0).
\end{eqnarray}
The action of these generators on the basis states $\{w_{m}^{j}\;|\; j=(0,
{1\over 2},1, \cdots),  -j \leq m \leq j\} $ may be taken as
\begin{eqnarray}
&&J_{+}.w_{m}^{j}=(j-m)(j+m+1)w_{m+1}^{j},\nonumber \\
&& J_{-}. w_{m}^{j}= w_{m-1}^{j}\;\; \hbox{for}\;\;\; m\geq -j+1, 
\qquad\qquad  
J_{-} .w_{-j}^{j}=0,\nonumber \\
&& J_{0} . w_{m}^{j}=2m\; w_{m}^{j}.
\end{eqnarray}
For $j={7\over 2}$, the representation $(2.37)$ may be immediately 
constructed by using the map $(2.38)$.

It is interesting to note that, following the map (2.38) an induced 
cocommutative coproduct structure (${\tilde \bigtriangleup}$) may be 
ascribed to the generators ($X,Y,H$) as follow
 \begin{eqnarray}   
&&{\tilde \bigtriangleup}( X) = {2\over h}\;\hbox{arctanh}({h\;
\bigtriangleup(J_{+})\over 2}), \nonumber\\
&& {\tilde \bigtriangleup}(Y)=\sqrt{1-{h^2\bigtriangleup(J_{+})^2 \over 4}}\;
\bigtriangleup(J_{-})\; \sqrt{1-{h^2\bigtriangleup(J_{+})^2 \over 4}}, 
\nonumber\\
&& {\tilde \bigtriangleup}(H)=\;\bigtriangleup(J_{0}).
\end{eqnarray}           
It is to be emphasized, however, that the coproduct structure 
${\tilde \bigtriangleup}$ treats the generators $(X,Y,H$) as elements of 
undeformed ${\cal U}(sl(2))$ algebra and is unrelated to the coproduct
structure of the Hopf algebra ${\cal U}_{h}(sl(2))$. Here we repeat the 
comment about usefulness of our systematic of development of the recursive 
scheme of finding the irreducible representations of the ${\cal U}_h(sl(2))$. 
This procedure is developed from the first principles and may be useful for
other nonlinear algebras, where the maps, similar to $(2.38)$, to the 
corresponding linear algebras may not be available.          

In the rest of the present section, we consider the ${\cal U}_h(e(2))$ 
algebra, which may be obtained \cite{Ball2} as a contraction of the 
 ${\cal U}_h(sl(2))$ algebra. Starting with an undeformed $e(2)$ generated by 
$(J, {\cal P}_{\pm})$ 
\begin{eqnarray}
&& [J , {\cal P}_{\pm} ] =\pm {\cal P}_{\pm}, \qquad \qquad
   [{\cal P}_{+} , {\cal P}_{-} ] =0,
\end{eqnarray}  
we, parallel to our previous prescription in $(2.38)$, define 
\begin{eqnarray}
&& \chi  = {2\over h} \hbox{arctanh} ({h {\cal P}_{+} \over 2}),\nonumber  \\
&& \eta = \sqrt{1- ({h\over 2} {\cal P}_{+})^{2}} \;\;{\cal P}_{-}\;\; 
  \sqrt{1- ({h\over 2} {\cal P}_{+})^{2}}= \biggl(1- ({h\over 2} 
{\cal P}_{+})^{2}\biggr)
\;{\cal P}_{-} ,\nonumber  \\
&& \zeta =2J.
\end{eqnarray}   
Then it follows
\begin{eqnarray}
&& [\zeta , \chi] = {2\over h} \sinh (h\chi), \qquad
 [\zeta,  \eta ] =-2 \;\eta \cosh (h\chi), \qquad
 [ \chi , \eta ]=0.
\end{eqnarray}
The Casimir is now 
\begin{eqnarray}
C = {\cal P}_{+}{\cal P}_{-} = {\eta \over h} \sinh  ({h\chi}). 
\end{eqnarray}
 This corresponds to the algebra considered in  \cite{Ball2}, except that we 
have not 
implicitly changed the signature of the metric from Euclidean to Lorentzian.   
The representations of $e(2)$ again, as before, provide explicitly those 
for the $h$-deformed case ${\cal U}_{h}(e(2))$. The standard unitary 
representations of $e(2)$ are, of course, infinite dimensional.

\sect{On ${\cal U}_h(e(3))$ and its Nonlinear Map}

One of the techniques used \cite{Cel1,Cel2} in the studies of the 
nonsemisimple 
quantized algebras is contraction. A singular transformation is performed in 
the vector space of the universal enveloping algebra of a suitable semisimple 
quantum algebra. These singular transformations were applied \cite{Cel1,Cel2}
to study various $q$-deformed inhomogeneous algebras. Another important 
application of this technique was made in the study of the $\kappa$-deformed
Poincar\'e algebra \cite{Luk1,Luk2}. Contractions can be implemented in 
different ways implying different consequences. There are two distinct classes.
In the first, $q$ goes to unity at the contraction limit, but, nonetheless, the
deformation persists, as the ratio of $(\ln q)$ with the contraction parameter 
remains finite. The above works \cite{Cel1,Cel2,Luk1,Luk2} are in this 
category. In an alternate scheme, contraction is performed \cite{Ac1,Ac2,Ac3},
 by retaining the value of $q$ (generic or root 
of unity) invariant. When $q$ is a root of unity, the second procedure is 
useful to construct, for example, the periodic representations for the 
contracted algebra \cite{Ac2}. In the present work, we limit our consideration 
to the first type. 

In the context the Jordanian 
deformations, a
construction ${\cal U}_h(so(4))={\cal U}_h(sl(2))\oplus {\cal U}_{-h}(sl(2))$
was used \cite{Ball} to obtain the quantized algebra ${\cal U}_h(so(4))$. The
$h$-deformed $3$-dimensional Euclidean algebra ${\cal U}_h(e(3))$, among 
others, was realized \cite{Ball,Shar} as a contraction of ${\cal U}_h(so(4))$
algebra. In \cite{Fad}, the connection between the classical Euclidean algebra 
$e(3)$ and the non-linear $\sigma$-model was pointed out. In that respect, the 
study of the various deformations of $e(3)$ assumes importance \cite{Cel2}. 
In searching for an alternate $h$-deformation of the $e(3)$ algebra, here we 
follow a closely 
parallel approach. Starting with a different choice of the generators of the
${\cal U}_h(so(4))$ algebra, we obtain a new ${\cal U}_h(e(3))$ algebra at
a contraction limit. The present ${\cal U}_h(e(3))$ has several interesting 
properties. In particular, the rotation algebra is preserved after the 
contraction. More importantly, reminiscent of the scenario discussed 
\cite{ACC} for the $U_{h}(sl(2))$ algebra, the generators of this 
${\cal U}_h(e(3))$ algebra may be realized, via a nonlinear map, in terms 
of the undeformed $e(3)$ generators. This tremendously simplifies the study 
of the irreducible representations of the proposed ${\cal U}_h(e(3))$ algebra.
Along similar lines, the corresponding $q$-deformed algebra ${\cal U}_q(e(3))$
was studied \cite{Ac4} before. The comparison between the two cases turns out 
to be of interest.

The $h$-deformed algebra ${\cal U}_h(so(4))$ is considered as a direct sum 
${\cal U}_h(so(4))={\cal U}_h(sl(2))\oplus {\cal U}_{-h}(sl(2))$. The choice of
the oppositely signed deformation parameters ($\pm h$) is necessary \cite{Cel2}
to avoid singularities in the coproducts after the contraction is achieved.
Let $(X_1,Y_1,H_1)$ and $(X_2,Y_2,H_2$) be two mutually commuting sets of 
generators, where each triplet satisfies the commutation relations(1.1). 
Their coalgebraic structure reads
\begin{eqnarray}
&& \bigtriangleup (X_{i})=X_i \otimes 1 + 1 \otimes X_i,\nonumber  \\
&& \bigtriangleup (Y_{i})=Y_i \otimes e^{h\theta_i X_i} + e^{-h\theta_i X_i} 
 \otimes Y_i, \nonumber\\
&& \bigtriangleup (H_{i})=H_i \otimes e^{h\theta_i X_i} + e^{-h\theta_i X_i} 
 \otimes H_i, \nonumber\\
&& \varepsilon (X_i)= \varepsilon (Y_i)= \varepsilon (H_i)=0,\nonumber \\
&& S(X_i) =- X_i,\nonumber \\
&&  S(Y_i) =- e^{h\theta_i X_i}Y_ie^{-h\theta_i X_i},\nonumber \\
&&  S(H_i) =- e^{h\theta_i X_i}H_ie^{-h\theta_i X_i},  
\end{eqnarray}                 
where $i=(1,2)$ and $\theta_1(\theta_2)=1(-1)$. In contrast to the usual 
practice \cite{Ball,Shar}, we make the following choice of the 
${\cal U}_h(so(4))$ generators
\begin{equation}
\begin{array}{ll}
 J_{+}=X_1+X_2, &\qquad\;\qquad  K_+=X_1-X_2,\nonumber \\
 J_- = Y_1e^{h X_2} + e^{-h X_1}Y_2, & \qquad \;\qquad
K_- =  Y_1e^{h X_2} - e^{-h X_1}Y_2,\nonumber \\
J_0 =H_1e^{h X_2} + e^{-h X_1}H_2,  & \qquad\qquad  \
K_0 =H_1e^{h X_2} - e^{-h X_1}H_2.
\end{array}
\end{equation}      
This choice is motivated by our intention to preserve a subalgebra $(J_{\pm},
J_0)$ satisfying the commutation relations of the ${\cal U}_h(sl(2))$
generators. The generators $(3.2)$ may be expressed in terms of the generators used in \cite{Ball,Shar}. The algebraic structure of ${\cal U}_h(so(4))$ now reads
\begin{eqnarray}
&& [J_0,J_+] = [K_{0},K_{+}]= {2\over h} \sinh (hJ_{+}),\nonumber  \\
&&  [J_0,J_-] = -J_{-} \cosh (hJ_{+})- \cosh (hJ_{+})\;J_{-},\nonumber  \\
&& [J_+,J_-] = [K_{+},K_{-}]=J_{0},\nonumber  \\
&&  [K_0,K_-] = -J_{-} e^{-hK_{+}}- e^{-hK_{+}} J_{-}- 
K_- \sinh (hJ_{+})- \sinh (hJ_{+})\; K_{-}, \nonumber \\
&& [J_0,K_+] =[K_0,J_+] = {2\over h} (\cosh(hJ_{+}) -e^{-hK_{+}}),\nonumber  \\
&& [J_0,K_-] = -K_{-} \cosh(hJ_{+}) -  \cosh(hJ_{+})\; K_{-}
-{h\over 8} \biggl((J_{0} +K_{0})e^{-h J_{+}}  \nonumber\\
&& \phantom{[J_0,K_-] = }+ e^{-h J_{+}} (J_{0} +
K_{0})\biggr) \biggl((J_{0} -K_{0})e^{h J_{+}} + e^{h J_{+}} (J_{0} -
K_{0})\biggr),\nonumber  \\
&& \nonumber \\
&& [K_0,J_-] = -K_{-} e^{-hK_{+}} -  e^{-hK_{+}} K_{-}
- \sinh (hJ_{+})\; J_{-} -J_{-}\; \sinh (hJ_{+}) + \nonumber  \\
&&  \phantom{[J_0,K_-] = }
{h\over 8} \biggl((J_{0} +K_{0})e^{-h J_{+}} + e^{-h J_{+}} (J_{0} +
K_{0})\biggr) \biggl((J_{0} -K_{0})e^{h J_{+}} + e^{h J_{+}} (J_{0} -
K_{0})\biggr), \nonumber  \\
&& [J_+,K_-] =[K_+,J_-] = K_{0},  \nonumber \\
&& [J_{+}, K_{+}] =0, \nonumber \\
&&  [J_-,K_-]= -{h \over 4} (J_- + K_-) 
\biggl( e^{-h J_{+}}( J_0 -K_0)  e^{h J_{+}}+
( J_0 -K_0)\biggr)  \nonumber \\
&&\phantom{ [J_-,K_-]= }    -{h \over 4} \biggl((J_0 + K_0) e^{-h J_{+}}
+  e^{-h J_{+}} (J_{0} +K_{0}) \biggr) (J_{-} -K_{-}) e^{h J_{+}} \nonumber \\
&& [J_{0}, K_{0}]= 2J_{0}\; \sinh (h J_{+}) +2K_{0} 
(e^{-h K_{+}} -\cosh (hJ_{+})).   
\end{eqnarray}    
One advantage of the algebra (3.3) in that the representation of the 
subalgebra consisting of the generators $(J_{\pm}, J_{3})$ is completely
known. This may help in finding the representation of the complete 
structure (3.3). The coalgebraic structure is as follows
\begin{eqnarray}
&& \bigtriangleup (J_{+})=J_{+} \otimes 1 + 1 \otimes J_{+},\nonumber   \\
&& \bigtriangleup (J_{-})=J_{-} \otimes \cosh (hJ_{+}) + 
e^{ -hK_{+}}\otimes J_{-} +K_{-} \otimes \sinh (hJ_{+}),\nonumber   \\
&& \bigtriangleup (J_{0})=J_{0} \otimes \cosh (hJ_{+}) + 
e^{ -hK_{+}}\otimes J_{0} +K_{0} \otimes \sinh (hJ_{+})   , \nonumber  \\
&& \bigtriangleup (K_{+})=K_{+} \otimes 1 + 1 \otimes K_{+}, \nonumber  \\
&& \bigtriangleup (K_{-})=K_{-} \otimes \cosh (hJ_{+}) + 
e^{ -hK_{+}}\otimes K_{-} +J_{-} \otimes \sinh (hJ_{+})   ,\nonumber   \\
&& \bigtriangleup (K_{0})=K_{0} \otimes \cosh (hJ_{+}) + 
e^{ -hK_{+}}\otimes K_{0} +J_{0} \otimes \sinh (hJ_{+})   , \nonumber  \\
&& \varepsilon (\xi)= 0 \qquad \hbox{for} \qquad \xi=(J_{+},
J_{-}J_{0},K_{+},K_{-},K_{0}),\nonumber   \\
&& S(J_{+}) =- J_{+},\nonumber   \\
&&  S(J_{-}) =- e^{hK_{+}}\biggl(J_{-} \cosh(hJ_{+}) - K_{-} \sinh (hJ_{+})
\biggr),\nonumber  \\
&& S(J_{0}) =- e^{hK_{+}}\biggl(J_{0} \cosh(hJ_{+}) - K_{0} \sinh (hJ_{+})
\biggr),
\nonumber  \\
&& S(K_{+}) =- K_{+},\nonumber   \\
&& S(K_{-}) =- e^{hK_{+}}\biggl(K_{-} \cosh(hJ_{+}) - J_{-} \sinh (hJ_{+})
\biggr),\nonumber  \\
&& S(K_{0}) =- e^{hK_{+}}\biggl(K_{0} \cosh(hJ_{+}) - J_{0} \sinh (hJ_{+})
\biggr).
\end{eqnarray}

The universal ${\cal R}$-matrix for the algebra ${\cal U}_h(sl(2))$ was
obtained \cite{Vlad,Shar} recently. For the algebra ${\cal U}_h(so(4))$,
the corresponding universal ${\cal R}$-matrix was constructed \cite{Shar} by
considering two copies of the ${\cal R}$-matrix of the algebra 
${\cal U}_h(sl(2))$. For our choice of the generators (3.2), the universal 
${\cal R}$-matrix for the ${\cal U}_h(so(4))$ algebra assumes the form 
\begin{eqnarray}  
&& {\cal R}=\exp\biggl[{h\over 8} (\Delta -\Delta') \nonumber  \\
&& \phantom{{\cal R}=}\biggl[{ 
(J_{0} + K_{0}) (1 - e^{-h (J_{+}-K_{+})}) (J_{+}+K_{+})
+  (J_{0}-K_{0}) (e^{h (J_{+}+K_{+})}-1)(J_{+}-K_{+})
\over \cosh (hJ_{+})- \cosh (hK_{+})} \biggr] 
\biggr],\nonumber \\
&& 
\end{eqnarray}  
where $\Delta'$ is the flipped coproduct map.

In order to construct the ${\cal U}_h(e(3))$ algebra by using the contraction 
procedure, we define a complex parameter $\omega$ and the generators 
$(P_{\pm},P_0)$ as follows
\begin{eqnarray}
&& \omega ={h\over \epsilon},\qquad \qquad P_{\pm} = \epsilon\;K_{\pm},
\qquad\qquad P_0=\epsilon\;K_0, 
\end{eqnarray}
while the generators $(J_{\pm},J_0)$ are not transformed. In the limit 
$\epsilon \longrightarrow 0$, the following algebraic structure of 
${\cal U}_h(e(3))$ is obtained from (3.3):
\begin{eqnarray}
&& [J_0, J_{\pm}] =\pm 2 J_{\pm},\qquad \qquad\qquad \qquad \qquad \;  
 [J_{+},J_{-}] =J_{0},\nonumber   \\
&&[P_{0}, P_{\pm}]=0,  \qquad \qquad \qquad \qquad \qquad \qquad [P_{+},P_{-}]=0,\nonumber \\
&&[J_{0},P_{+}]=[P_{0},J_{+}]={2\over \omega} (1-e^{-\omega P_{+}}), 
\nonumber \\
&& [J_{0}, P_{-}] = -2P_{-} +{\omega \over 2} P_{0}^{2},  \nonumber \\
&&[P_{0}, J_{-}] = -2e^{-\omega P_{+}} P_{-} -{\omega \over 2} P_{0}^{2},
\nonumber \\ 
&& [J_{+}, P_{-}] = [P_{+}, J_{-}]=P_{0}, \nonumber  \\
&&[J_{+}, P_{+}] = 0, \qquad \qquad \qquad \qquad \qquad \qquad [J_{-}, P_{-}] = \omega P_{0}P_{-}, \nonumber\\
 &&[J_{0}, P_{0}] =-2 P_{0} (1-e^{-\omega P_{+}}).  
\end{eqnarray}
A distinctive feature of this algebra is that, unlike the results obtained 
in \cite{Shar} for the corresponding case, here the rotation group is 
preserved. The coalgebra maps for ${\cal U}_h(e(3))$ are read from (3.4) 
at the contraction limit
\begin{eqnarray}
&& \bigtriangleup (J_{+})=J_{+} \otimes 1 + 1 \otimes J_{+}, \nonumber  \\
&& \bigtriangleup (J_{-})=J_{-} \otimes 1 + 
e^{ -\omega P_{+}}\otimes J_{-} +\omega P_{-} \otimes J_{+}   ,  \nonumber \\
&& \bigtriangleup (J_{0})=J_{0} \otimes 1 + 
e^{ -\omega P_{+}}\otimes J_{0} +\omega P_{0} \otimes J_{+}   , \nonumber  \\
&& \bigtriangleup (P_{+})=P_{+} \otimes 1 + 1 \otimes P_{+}, \nonumber  \\
&& \bigtriangleup (P_{-})=P_{-} \otimes 1 + 
e^{ -\omega P_{+}}\otimes P_{-}, \nonumber  \\
&& \bigtriangleup (P_{0})=P_{0} \otimes 1 + 
e^{ -\omega P_{+}}\otimes P_{0} ,  \nonumber \\
&& \varepsilon (\xi)= 0 \qquad \hbox{for} \qquad \xi=(J_{+},
J_{-}, J_{0},P_{+},P_{-},P_{0}), \nonumber  \\
&& S(J_{+}) =- J_{+}, \qquad 
S(J_{-}) =- e^{\omega P_{+}}(J_{-} - \omega P_{-} J_{+}), \qquad 
 S(J_{0}) =- e^{\omega P_{+}}(J_{0} - \omega P_{0} J_{+}), \nonumber \\
&& S(P_{+}) =- P_{+}, \qquad  
S(P_{-}) =- e^{\omega P_{+}}P_{-}, \qquad \qquad\;\;\;\qquad 
 S(P_{0}) =- e^{\omega P_{+}} P_{0}. 
\end{eqnarray}   
The universal ${\cal R}$-matrix for the ${\cal U}_h(e(3))$ algebra may be 
obtained from $(3.5)$ at the $\epsilon \longrightarrow 0$ limit as the 
singular terms cancel:
\begin{eqnarray}
{\cal R} = \exp\biggl[ {\omega \over 4} (\bigtriangleup - \bigtriangleup')
\biggl[ {(P_{0} J_{+} + J_{0} P_{+}) (1-e^{\omega P_{+}}) + \omega P_{0} 
P_{+} J_{+} e^{\omega P_{+}} \over 1 - \cosh ( \omega P_{+})}\biggr]\biggr].
\end{eqnarray}   

Closely paralleling our earlier description of the ${\cal U}_h(sl(2))$ algebra,
we demonstrate here that the algebra (3.7) may be mapped on the classical 
$e(3)$ algebra. To this end, we define the generators 
\begin{eqnarray}
\Pi_{+} = {1\over \omega} (e^{\omega P_{+}}-1 ), \qquad 
\Pi_{-} = P_{-} -{\omega \over 4} P_{0}^{2} e^{\omega P_{+}}, \qquad  
\Pi_{0} = P_{0} e^{\omega P_{+}}.
\end{eqnarray}   
From the algebra (3.7), it follows that ($J_{\pm},J_{0}, \Pi_{\pm}, \Pi_{0}$)
obey classical algebra $e(3)$ algebra, where $(\Pi_{\pm}, \Pi_{0})$ play the 
role of generators of translations. We obtain
\begin{eqnarray}
&& [J_{0},J_{\pm}] =\pm 2 J_{\pm},  \qquad\qquad [J_{+}, J_{-}]= J_{0}, 
\qquad\qquad [\Pi_{0}, \Pi_{\pm} ] =0, \qquad\qquad 
[\Pi_{+},\Pi_{-}]=0,\nonumber\\  
&& [J_{+}, \Pi_{+}]=0,   \;\;\;\;\;\;\qquad \qquad [J_{+}, \Pi_{-}] =\Pi_{0}, 
\;\;\;\;\;\qquad [J_{+}, \Pi_{0} ] = -2 \Pi_{+},  \nonumber \\
&& [J_{-}, \Pi_{+}] = - \Pi_{0}, \qquad\; \qquad
[J_{-}, \Pi_{-}] =0, \qquad \qquad [J_{-}, \Pi_{0}] = 2 \Pi_{-}, 
 \nonumber \\
&& [J_{0}, \Pi_{\pm}] = \pm 2 \Pi_{\pm}, \qquad \qquad  [J_{0},\Pi_{0}]=0,     
\end{eqnarray} 
The Casimir operators are
\begin{eqnarray}
&& C_1 = \Pi_{+} \Pi_{-}+{1\over 4} \Pi_{0}^{2}, \qquad \qquad \nonumber\\
&& \phantom{C_1}= {1 \over \omega}(e^{\omega P_{+}} -1)
(P_{-} - { \omega \over 4} P_{0}^{2}  e^{\omega P_{+}})+ 
{ 1 \over 4} P_{0}^{2}  e^{2\omega P_{+}}, \\ 
&& C_2 = J_{+}\Pi_{-} + J_{-} \Pi_{+} + {1\over 2} J_{0}\Pi_{0} \nonumber \\
&&  \phantom{C_1}=J_{+} (P_{-}-{\omega\over 4} P_{0}^{2} e^{\omega P_{+}})
+  {1 \over \omega}J_{-}(e^{\omega P_{+}} -1) + {1 \over 2} J_{0}P_{0}
e^{\omega P_{+}}.       
\end{eqnarray}
The representation theory of the classical $e(3)$ algebra may now be 
readily used to obtain the representations of the algebra (3.7). For the
classical $e(3)$ generated by $(J_{\pm},J_{0}, \Pi_{\pm}, \Pi_{0}$) we
may introduce the standard representations in the momentum or the angular 
momentum basis. Concerning the action of the generators $(P_{\pm}, P_{3})$
on momentum bases, we note the following. The inverse map of (3.10) reads 
\begin{eqnarray}     
P_{+} = {1\over \omega} \ln (1+\omega \Pi_{+}), \qquad 
P_{-} = \Pi_{-} + {\omega \over 4} \;{\Pi^{2}_{0} \over (1 +\omega \Pi_{+})}, 
\qquad P_{0} = {\Pi_{0} \over  (1 +\omega \Pi_{+}) }
\end{eqnarray}
On the momentum basis (where the eigenvalues of $\Pi_{\pm}$ and $\Pi_{0}$ are
taken to be ${\tilde \Pi}_{\pm}$ and ${\tilde \Pi}_{0}$  respectively), it is 
particularly evident from $(3.14)$ that the eigenvalues of the operators 
$(P_{\pm}, P_{0})$ diverge for $(1+\omega{\tilde \Pi}_{+})=0$. Moreover the 
eigenvalue of $P_{+}$ develops an imaginary part for $(1+\omega{\tilde 
\Pi_{+}})< 0$. In fact the
whole real domain for the eigenvalues of operators $(P_{\pm}, P_{0})$ 
correspond to the restricted domain of the classical eigenvalues given by 
$(1 +\omega {\tilde \Pi}_{+}) > 0$. 

The singularities in the inverse map $(3.14)$ contrast very sharply with a
similar construction obtained by one of us \cite{Ac4} in the context of 
$q$-deformed algebra ${\cal U}_q(e(3))$. We will briefly review the mapping 
of the nonlinear algebraic structure of ${\cal U}_q(e(3))$ on the classical 
$e(3)$ algebra in the Appendix. From $(A.4)$, it is evident that, there the 
nonlinearity enters in the map through the positive definite `invariant 
mass' $C_1$, and, consequently, the map is invertible without any 
singularity. It is known \cite{Kosi} that the $\kappa$-deformed Poincar\'e 
algebra may be realized in terms of the generators of the classical 
Poincar\'e algebra. We note the analogies and contrasts between the map
\cite{Kosi} for the $\kappa$-Poincar\'e algebra and the corresponding maps 
for the algebras ${\cal U}_h(e(3))$ and ${\cal U}_q(e(3))$. As noted earlier, 
the nonlinearity enters the map $(3.14)$ for the ${\cal U}_h(e(3))$ algebra
through the classical momentum operator, which is not positive definite, and, 
therefore, the map shows singularity. As for the ${\cal U}_q(e(3))$ example 
discussed in the Appendix, there is no singularity in the map. For the 
$\kappa$ deformed algebra, the nonlinearity enters \cite{Kosi} 
through the classical energy operator and, for the positive 
energy solutions, the map shows no singularity.

\sect{Conclusion}

We obtained the finite dimensional highest weight representations of 
${\cal U}_h(sl(2))$ algebra as factor representations by using the standard 
singular vector treatment. This was done by solving a set of recursion 
relations valid for the matrix elements of the ${\cal U}_h(sl(2))$ generators.
These irreducible representations may also be determined by mapping the 
${\cal U}_h(sl(2))$ algebra on the classical $sl(2)$ algebra. The 
corresponding map of the ${\cal U}_h(e(2))$ algebra was obtained by using
a contraction method. Taking a nonlinear combination of generators of two 
copies of a $h$-deformed $sl(2)$ algebra, the ${\cal U}_h(so(4))$ algebra 
was constructed. An ${\cal U}_h(e(3))$ algebra was constructed by contracting 
${\cal U}_h(so(4))$. The ${\cal U}_h(e(3))$ may be mapped on the classical 
$e(3)$ algebra. A physically interesting feature is that this map, unlike the
previously known cases of ${\cal U}_q(e(3))$ algebra and the 
$\kappa$-Poincar\'e algebra, exhibits a singular behavior. 
In the momentum basis, a restricted domain of the eigenvalues of the 
classical operators is mapped on the whole real domain of the eigenvalues of
the deformed operators.

\vskip 2cm 
\noindent {\bf Acknowledgments:}

One of us (RC) wants to thank A. Chakrabarti for a kind invitation. He is 
also grateful to the members of the CPTH group for their kind hospitality.

\newpage

\appendix
\def\thesection{Appendix}

\section{}

Here we deviate from the main body of the paper and give a summary of the 
construction of the $q$-deformed algebra ${\cal U}_q(e(3))$, where the 
algebraic structure may also be mapped \cite{Ac4} on the undeformed $e(3)$ 
algebra. Unlike the $h$-deformed case discussed in $(3.10)$ and $(3.14)$, 
the map obtained in \cite{Ac4} is invertible and non-singular.

Starting with a contraction ${\cal U}_q(so(4))= {\cal U}_q(sl(2))\oplus 
{\cal U}_{-q}(sl(2))$ where $({\hat J}_ {\pm}^{(i)}, {\hat J}_ {0}^{(i)})$ are the generators of 
the two deformed $sl(2)$ algebras, we define 
\begin{eqnarray}
&& {\hat J}_ {0}= {\hat J}_ {0}^{(1)}+ {\hat J}_ {0}^{(2)}, \qquad   
{\hat J}_ {\pm}= {\hat J}_ {\pm}^{(1)}q^{-{\hat J}_ {0}^{(2)}}+ 
{\hat J}_ {0}^{(2)}q^{{\hat J}_ {0}^{(1)}}, \nonumber \\
&& {\hat K}_{0}= {\hat J}_ {0}^{(1)}- {\hat J}_ {0}^{(2)}, \qquad   
{\hat K}_{\pm}= {\hat J}_ {\pm}^{(1)}q^{-{\hat J}_ {0}^{(2)}}- 
{\hat J}_ {0}^{(2)}q^{{\hat J}_ {0}^{(1)}}.
\end{eqnarray}        
The modifications with respect to \cite{Cel2}, in the definitions of
${\hat J_{\pm}}$, ${\hat K_{\pm}}$ are suitable for our purpose.   
We now perform the following contraction 
\begin{eqnarray}
\displaystyle {\hat K}_{\pm} = {{\hat P}_{\pm} \over \epsilon}, \qquad  
{\hat K}_{0} = {{\hat P}_{0} \over \epsilon},\qquad q=e^{\epsilon \Omega} 
\end{eqnarray}
and take the limit $\epsilon \longrightarrow 0$. The  ${\cal U}_q(so(4))$ 
algebra constructed with the generators $(A.1)$ preserves the rotation 
subalgebra, both before and after the contraction. It also permits to exhibit 
a nonlinear mapping to the classical algebra in a direct fashion. After 
contraction, the ${\cal U}_q(e(3))$ algebra reads
\begin{eqnarray}     
&& [{\hat J}_ {0}, {\hat J}_ {\pm}] =\pm 2 {\hat J}_ {\pm}, 
\qquad\qquad\;\;\;\qquad [{\hat J}_ {+}, {\hat J}_ {-}]= {\hat J}_ {0}, \nonumber\\
&& [{\hat P}_{0}, {\hat P}_{\pm} ] =0, \qquad\qquad\;\;\qquad \qquad 
[{\hat P}_{+},{\hat P}_{-}]=0, \nonumber \\
&&[{\hat J}_{0},{\hat P}_{\pm}]= [{\hat P}_{0}, 
{\hat J}_{\pm}]=\pm 2 {\hat P}_{\pm}, \qquad
  [{\hat J}_{0},{\hat P}_{0}]=0,\nonumber\\
&& [{\hat P}_{\pm}, {\hat J}_ {\pm}]=\pm \Omega {\hat P}_{\pm}^{2}
,\qquad \qquad\qquad [{\hat J}_ {\pm}, {\hat P}_{\mp}] =\pm {1\over \Omega} 
(e^{\Omega {\hat P}_{0}} -1). 
\end{eqnarray}  
The algebra $(A.3)$ may be mapped on the classical $e(3)$ algebra. Let 
$({\hat J}_{\pm}, {\hat J}_{0}, {\hat \Pi}_{\pm},{\hat \Pi}_{0})$ be the 
generators of the classical $e(3)$ algebra. Then the map reads 
\begin{eqnarray}
&& e^{-{\Omega \over 2}{\hat P}_{0}} = (1+ {1\over 2} C_1 \Omega^{2}) -
{\Omega \over 2} (1+{1\over 4}C_1\Omega^{2})^{1/2}{\hat\Pi}_{0},\nonumber\\
&& {\hat P}_{\pm}e^{-{\Omega \over 2}{\hat P}_{0}}=(1+{1\over 4}C_1 
\Omega^{2})^{1/2} {\hat\Pi}_{\pm}  ,
\end{eqnarray}  
where $C_1$ is the positive definite `invariant mass' of the undeformed $e(3)$ 
algebra. In terms 
of the generators of ${\cal U}_q(e(3))$, the Casimir operator $C_1$  
($={\hat \Pi}_+{\hat \Pi_-}
+{1\over 4}{\hat \Pi_{0}}^2$) has the form 
\begin{eqnarray}
&& C_1 = {\hat P}_{+} {\hat P}_{-} e^{-{\Omega \over 2}{\hat P}_{0}} +{1\over \Omega^{2}}
(e^{{\Omega \over 2}{\hat P}_{0}}+ e^{-{\Omega \over 2}{\hat P}_{0}}-2).
\end{eqnarray}     
Therefore the mapping (A.3) is invertible without any singularity. This 
clearly contrasts the map described in Section 3 for the $h$-deformed algebra
${\cal U}_h(e(3))$.

\end{document}